\documentclass[11pt,a4paper]{amsart}
\usepackage[english]{babel}
\usepackage[utf8]{inputenc} 
\usepackage[T1]{fontenc}
\usepackage{verbatim}
\usepackage{palatino}
\usepackage{amsmath}
\usepackage{mathabx}
\usepackage{amsbsy}
\usepackage{amssymb}
\usepackage{amsthm}
\usepackage{amsfonts}
\usepackage{graphicx}
\usepackage{esint}
\usepackage{color}
\usepackage{mathtools}
\usepackage{overpic}
\usepackage{mathrsfs}
\usepackage{tabularx} 
\mathtoolsset{showonlyrefs}

\usepackage[colorlinks = true, citecolor = black]{hyperref}
\author{Andr\'as M\'ath\'e \& William O'Regan}
\title[Discretised sum-product theorems by Shannon-type inequalities]{Discretised sum-product theorems\\by Shannon-type inequalities}
\address[$\ddag$]{Zeeman Building,
University of Warwick,
Coventry,
West Midlands,
United Kingdom}
\email{a.mathe@warwick.ac.uk}
\address[$\dagger$]{1984 Mathematics Road, Vancouver, British Columbia, V6T 1Z2, Canada}
\email{woregan@math.ubc.ca}
\date{\today}
\subjclass[2010]{05B99, 28A78, 28A80}
\keywords{Discretised sum-product, discretised ring conjecture, Shannon entropy, Pl\"unnecke--Ruzsa, projection theorems}
\thanks{A.M. is supported by the Hungarian National Research, Development and Innovation Office -- NKFIH\-, 124749. This work was completed while W.O.R. was supported by the EPSRC via the project \emph{Ergodic and combinatorial methods in fractal geometry}, project ref.~2443767.}

\newcommand{\E}{\mathrm{H}}
\newcommand{\Ed}{\E_\delta}

\newcommand{\D}{\mathcal{D}}

\newcommand{\om}{\omega}

\newcommand{\R}{\mathbb{R}}
\newcommand{\pr}{\mathbb{P}}

\newcommand{\N}{\mathbb{N}}

\newcommand{\Z}{\mathbb{Z}}

\newcommand{\calL}{\mathcal{L}}

\newcommand{\calD}{\mathcal{D}}

\newcommand{\lbox}{\underline\dim_{\mathrm{B}}}
\newcommand{\ubox}{\overline\dim_{\mathrm{B}}}

\newcommand{\spt}{\operatorname{spt}}

\newcommand{\dimh}{\dim_{\rm{H}}}

\newcommand{\dist}{\operatorname{dist}}

\newcommand{\Om}{\Omega}

\numberwithin{equation}{section}
\theoremstyle{plain}
\newtheorem{thm}[equation]{Theorem}

\newtheorem{lemma}[equation]{Lemma}

\newtheorem{prop}[equation]{Proposition}

\theoremstyle{definition}

\newtheorem{definition}[equation]{Definition}

\theoremstyle{remark}
\newtheorem{remark}[equation]{Remark}

\newlength\tindent

\renewcommand{\indent}{\hspace*{\tindent}}
\newcommand{\nref}[1]{(\hyperref[#1]{#1})}
\begin{document}
\begin{abstract}
By making use of arithmetic information inequalities, we give a strong quantitative bound for the discretised ring theorem. In particular, we show that if $A \subset [1,2]$ is a $(\delta,\sigma)$-set, with $|A| = \delta^{-\sigma},$ then $A+A$ or $AA$ has $\delta$-covering number at least $\delta^{-c}|A|$ for any $0 <  c < \min\{\sigma/6, (1-\sigma)/6\}$ provided that $\delta > 0$ is small enough.
\end{abstract}
\maketitle
\section{Introduction}
Erd\H{o}s and Volkmann in \cite{erd} showed that for any $\sigma\in [0, 1]$ there exists a Borel subgroup of the reals  with Hausdorff dimension $\sigma.$ They conjectured that the same does not hold for Borel subrings, more, there does not exist a Borel subring of the reals with Hausdorff dimension strictly between zero and one. Their conjecture was proved by Edgar and Miller in \cite{ed}, using projection theorems of fractal sets. Essentially at the same time, Bourgain \cite{bou03} independently proved the conjecture via solving the discretised ring conjecture of Katz and Tao \cite{kat}. 

\indent A classical example of the occurrence of sum-product phenomena is the following theorem from  Erd\H{o}s and Szemer\'edi \cite{erdsz}. They state that there exists an $\epsilon > 0$ and a $C_\epsilon > 0$ such that for every finite subset of integers $A$ at least one of the sumset $A+A$ or the product set $AA$ is large in the sense that
\begin{displaymath}
\max(|A+A|, \ |AA|) \ge C_\epsilon |A|^{1+\epsilon}.
\end{displaymath}
Indeed, this asserts that any finite subset of the integers can not even approximately resemble the structure of a ring. They conjectured that a positive constant $C_\epsilon$ exists for every $0 < \epsilon<1,$ that is, at least one of $|A + A|$ or $|AA|$ must be nearly as large as possible.

\indent The discretised sum-product problem (or discretised ring problem) of Katz--Tao \cite{kat} is the discretised version of the fractal analogue of the Erd\H{o}s--Szerem\'edi problem. Vaguely, it asks/asserts that if $A\subset \R$ behaves like an $\sigma$-dimensional set at scale $\delta$ in a certain sense, then at least one of $A+A$ and $AA$ behaves like an $(\sigma+c)$-dimensional set at scale $\delta$ (in a  different and slightly weaker sense),  where the positive constant $c$ should depend only on $\sigma$. As previously mentioned, it was first proved in 2003 by Bourgain in \cite{bou03}, and represented again with weaker non-concentration conditions by Bourgain--Gamburd in 2008, \cite{bougam} and Bourgain in 2010 in \cite{bou}. No explicit bound on the constant was presented. Further examination of Bourgain's papers would suggest that the explicit constant gained following his exact method would be very small. Strong explicit constants were gained by Guth, Katz, and Zahl \cite{gut}, by Chen in \cite{che}, and Fu and Ren \cite[Corollary 1.7]{furen}. 

\indent The discretised sum-product also has many other applications.  For instance it is closely related to the Falconer distance set problem and the dimension of Furstenberg sets, see Katz and Tao \cite{kat} for more details.  For some applications of discretised sets to projections of fractal sets see, for example, He \cite{he}, Orponen \cite{orp2}, \cite{orpabcd}, and Orponen--Shmerkin--Wang \cite{orp} and the references therein. For the applications of discretised sum-product to the Fourier decay of measures see Li \cite{li}.  

\indent The aim of this paper is to provide a strong bound for $c$ for the Katz--Tao discretised sum-product problem. We show that $c$ can be taken arbitrarily close to $\sigma/6$ if $\sigma \leq 1/2$ and arbitrarily close to $(1-\sigma)/6$ when $1/2<\sigma<1$. 

\indent Clearly, $c$ cannot exceed $\sigma$ nor $1-\sigma$. It is unclear if it is reasonable to conjecture that $c$ can be taken to be (nearly) $\sigma$ when $\sigma$ is small (analogously to the Erd\H{o}s--Szemer\'edi conjecture). On the other hand, when $\sigma>1/3$, $c$ cannot be larger than $(1-\sigma)/2$, see Proposition \ref{prop.upperbound}.

\indent The approach in this paper is to start with theorems from fractal geometry that imply that certain arithmetic operations necessarily increase the dimension of any set $A\subset \R$ and then to use information inequalities to extract that simpler arithmetic operations (in this case, addition and multiplication) must already increase the dimension. Bourgain's original proof of the discretised ring conjecture and many improvements since relied on theorems of additive combinatorics (Ruzsa and Pl\"unnecke--Ruzsa inequalities). Our information inequalities make use of both the additive and multiplicative structure of the underlying field. All these inequalities are immediate corollaries of certain instances of the submodularity inequality, that is, that the conditional mutual information of two random variables given a third is non-negative.

\indent The Shannon entropy version of the Pl\"unnecke--Ruzsa inequalities were obtained by \newline \cite{mad}; see also \cite{taoent}. Our proof relies on a recent theorem of Orponen--Shmerkin--Wang in fractal geometry \cite{orp}. Their theorem is a generalisation of a classical theorem of Marstrand. See §\ref{proofsketch} below for details and a brief informal overview of our proofs.

\indent Since the first version of this preprint was uploaded, similar bounds were obtained in \cite{orpshabc} and then improved further in \cite{renwang}.
\subsection{Definitions and notation} The function $\log$ will always be to base 2. For some compact $A \subset \R^d$ we denote $\dimh A$ to be its Hausdorff dimension, $\lbox A, \ubox A $ to be its lower and upper box dimensions. For a measure $\mu$ on a space $X$ we define the \textit{conditioned measure} of $\mu$ with respect to $Y \subset X$ by $\mu_{|Y}(A) = \frac{\mu(A \cap Y)}{\mu(Y)}$ for all measurable $A \subset X.$ Let $C > 1$ and $0 < s < d.$ We say that a measure $\mu$ on $\R^d$ is $(s,C)$\textit{-Frostman} if it is a Radon probability measure with the non-concentration condition
\begin{equation}
    \mu(B(x,r)) \leq Cr^s \text{ for all } x \in \R^d, r > 0.
\end{equation}

\indent In the below $\delta >0$ is the \textit{scale} in which we will view our fractal set. The functions $f,g$ will be some quantity relating to the scale we are at, for example, $f 
= N_\delta(A),$ where $A$ is the set we are examining. For such an $f,$ we wish to understand for which $\alpha > 0$ we have that $f$ `behaves' like $\delta^{-\alpha}.$ The notation below makes this precise. Fix an exponent $0< \sigma <1,$ and a roughness parameter $C > 0.$ For two functions $f,g:(0,1] \rightarrow [0,\infty)$ we write $f \lesssim g$ if there exists a constant $K > 0,$ (it may depend on $C, d,$ and $\sigma$) so that 
\begin{equation}
    f(\delta) \leq Kg(\delta) \text{ for all scales }\delta \in (0,1].
\end{equation}
We write $f \gtrsim g$ if $g \lesssim f.$ We write $f \sim g$ if $f \gtrsim g$ and $g \gtrsim f.$

 \indent We write $f \lessapprox g$ if for all $0 < \epsilon < 1$ there exists a constant $K_\epsilon > 0$ (it may depend on $C,d, \epsilon, \sigma),$ so that
\begin{equation}
    f(\delta) \leq K_\epsilon \delta^{-\epsilon} g(\delta) \text{ for all scales } \delta \in (0,1].
\end{equation}
We write $f \gtrapprox g$ if $g \lessapprox f.$ We write $f \approx g$ if $f \gtrapprox g$ and $g \gtrapprox f.$
For example: $C \sim C^{\sigma} \sim 1;$ if $A \subset \R$ has box dimension $\sigma,$ then $N_\delta(A) \approx \delta^{-\sigma}.$ \textbf{Main point}: the implicit constants may \textbf{not} depend on the scale we are working with. 
\begin{definition}[$(\delta,\sigma, C)$-set]\label{deltasigma.def}
Fix a scale $\delta > 0,$ and a constant $C > 0.$ We say that a finite non-empty $\delta$-separated set $A \subset \R$ is a $(\delta,\sigma, C)$-\textit{set} if $A$ satisfies the following non-concentration condition:
     \begin{displaymath}
  |A \cap B(x,r)| \leq Cr^\sigma|A| \qquad x \in \R^d, \ r \geq \delta.   
\end{displaymath}
\end{definition}
 We remark that by setting $r = \delta$ that we must have $|A| \geq \delta^{-\sigma}/C.$ 
 
\indent A $(\delta, \sigma,C)$-set can be considered as the discrete approximation of a set in $\R$ with `dimension' $\sigma$ at scale $\delta.$ See \cite[Lemma 5.3]{bflm} for a precise formulation. 
\subsection{Main results}
In the below $\pm$ means the result is true for both $+$ and $-.$ For a bounded $A \subset \R,$ $N_\delta(A)$ denotes the least number of intervals of length $\delta$ needed to cover $A.$ The main results are the following. 
\begin{thm}\label{sumprod}
Let $0 < \delta, \sigma < 1, C >0.$ For all $(\delta, \sigma,C)$-sets $A \subset [1,2]$ we have
\begin{align}
N_\delta(A\pm A)^2N_\delta(A A)^4 &\gtrapprox \delta^{-5\sigma - c},\\
N_\delta(A\pm A)^2N_\delta(A/A)^3 &\gtrapprox \delta^{-4\sigma - c} \label{eq.hard2},
\end{align}
where $c = \min\{2\sigma,1\}.$ The implicit constants depend on $C$ and $\sigma.$
\end{thm} 
As a simple corollary we get the discretised ring theorem. 
\begin{thm}\label{disring}
Let $0 < \delta, \sigma < 1, C >0.$ For all $(\delta, \sigma,C)$-sets $A \subset [1,2]$ we have
\begin{align*}
N_\delta(A\pm A) + N_\delta(A A) &\gtrapprox \delta^{-\sigma - c},\\
N_\delta(A\pm A) + N_\delta(A/A) &\gtrapprox \delta^{-\sigma -c'},
\end{align*}
where $c = \min\{\sigma,1-\sigma\}/6$ and $c' = \min\{\sigma,1-\sigma\}/5.$
\end{thm}
Stated in terms of dimension we are also able to get the following. The proof follows from an application of \cite[Lemma 5.3]{bflm} along with basic properties of Hausdorff dimension. 
\begin{thm}\label{sumproddim}
For all $0 < s < 1$ and for all Borel sets $A \subset \R$ with Hausdorff dimension $s$ we have the following:
\begin{align*}
\lbox ((A A)^4 \times (A\pm A)^2) &\geq 5s+ \min\{2s,1\},\\
\lbox ((A/A)^3 \times (A\pm A)^2) &\geq 4s + \min\{2s,1\}.
\end{align*}
\end{thm}
Here $(AA)^4, (A/A)^3, (A+A)^2, (A-A)^2$ denote Cartesian products. 
Using product formulae the following follows from Theorem \ref{sumproddim} immediately. 
\begin{thm}
    For all $0 < s < 1$ and for all Borel sets $A \subset \R$ with Hausdorff dimension $s$ we have the following:
\begin{align*}
\max\{\lbox(A \pm A), \ubox(A A)\} &\geq \min\{7s/6, (5s+1)/6\},\\
\max\{\ubox(A \pm A), \lbox(A A)\} &\geq \min\{7s/6, (5s+1)/6\},\\
\max\{\lbox(A \pm A), \ubox(A/A)\} &\geq \min\{6s/5, (4s+1)/5\},\\
\max\{\ubox(A \pm A), \lbox(A/A)\} &\geq \min\{6s/5, (4s+1)/5\}.
\end{align*}
\end{thm}
\subsection{Proof sketch}\label{proofsketch}
We will rely on a recent generalisation of Marstrand's theorem by Orponen--Shmerkin--Wang \cite{orp}. They proved that for every pair of Borel sets $E,F \subset \R^2$ which are both not contained in a line, and both of Hausdorff dimension $s\in (0,2)$, the set of directions between $E$ and $F$ (that is, directions of line segments with one endpoint in $E$ and one endpoint in $F$) has Hausdorff dimension at least $\min\{1,s\}$ (and has positive Lebesgue measure if $s>1$). (We will need their stronger version of the same theorem involving Frostman estimates.)

\indent Now let $A\subset [1,2]$ be a Borel set of Hausdorff dimension $s.$ Let $E=A\times A$ and $F=(-A)\times (-A)$. Then both $E$ and $F$ have Hausdorff dimension at least $\min\{1,2s\}$ and the set of directions (and slopes) realised by line segments connecting a point of $E$ to a point of $F$ has Hausdorff dimension at least $\min\{1,2s\}$. Thus
\begin{equation}\label{ep}
\left\{\frac{a+b}{c+d}\in \R\,:\,a,b,c,d\in A\right\}
\end{equation}
has Hausdorff dimension at least $\min\{1,2s\}$. (This is noted in \cite{orp}.)

\indent Let $X,Y,Z,W$ be independent identically distributed random variables taking values in $A$ with an appropriate distribution (a Frostman measure on $A,$ for example). Then by submodularity, we have
$$\E\left(\frac{X+Y}{Z+W}\right) + 5\E(X) \le 2\E(X+Y) + 4\E(XY).$$
Here $\E$ is an appropriate version of Shannon entropy. In terms of `dimension', this inequality intuitively means that
$$\dim \left(\frac{A+A}{A+A}\right)+ 5\dim(A) \le 2\dim(A+A) + 4\dim(AA).$$
By \eqref{ep}
\begin{displaymath}
\min\{1,2s\}+5s\le 2\dim(A+A) + 4\dim(AA).
\end{displaymath}
In particular, at least one of $A+A$ and $AA$ should have `dimension' at least $\min\{(1+5s)/6,7s/6\}$. 

\begin{remark}
Our bounds for the sum--product problem depend on the arithmetic information inequalities we have found.
Given another information inequality, provided that one has a good lower bound for the left hand side (in the continuous/discretised setting as in the examples above), results for fractal sets and discretised $(\delta,\sigma)$-sets can be readily obtained by following similar approaches as presented in this paper.
\end{remark}
\begin{remark}\label{careful}
One has to be careful with stating the sum-product problem for fractal dimensions. In particular, the naive sum-product conjecture for Hausdorff dimension fails: for every $0<\sigma < 1$ there are compact sets $A\subset \R$ of Hausdorff dimension $\sigma$ such that both $A+A$ and $AA$ have Hausdorff dimension $\sigma$. The problem is that $A+A$ and $AA$ can be small at different scales, which is enough to make their Hausdorff dimension small. See Section \ref{section.examples}.
\end{remark} 
\subsection*{Acknowledgements}  Thanks are given to Tim Austin and Tam\'as Keleti for reading an earlier version of this manuscript. We thank the anonymous referee, whose comments greatly improved the quality of the exposition. 
\section{Preliminaries}
\subsection{Geometric measure theory}
We will need the following stability property, which is a straightforward property of $(\delta,\sigma,C)$-sets. We include the straightforward proof. 

\indent For two non-empty compact subsets $A,B \subset \R^2$ we define $\dist(A,B),$ a measure of separation between $A$ and $B,$ by 
\begin{equation}\label{eq.distsets}
\dist(A,B) = \inf\{d(a,b): a \in A, b \in B\}.
\end{equation}
\begin{lemma}\label{lem.twosubsets}
    Let $A \subset [1,2]$ be a $(\delta, \sigma,C)$-set. Set $C' = (6C)^{1/\sigma+1}.$ Then there exists $A_1, A_2 \subset A,$ with both $A_1$ and $A_2$ being $(\delta,\sigma,C')$-sets, and $\dist(A_1,A_2) \geq C'^{-1}.$
\end{lemma}
\begin{proof}
    Fix $\rho = (6C)^{-1/\sigma}.$ Let $\calD$ be the collection disjoint intervals of length $\rho$ intersecting $[1,2].$ By the pigeonhole principle, fix $I \in \D$ with $|A \cap I| > \rho^{-1}|A|.$ We have
    \begin{align*}
|A| &\leq \sum_{J \in \D}|A \cap J| \\
&= |A \cap I| + \sum_{J \text{ adjacent to } I}|A \cap J| + \sum_{J \text{ not adjacent to } I}|A \cap J| \\
&\leq 3C\rho^\sigma|A| +  \sum_{J \text{ not adjacent to } I}|A \cap J|.
    \end{align*}
Therefore 
\begin{equation}
    |A|/2 \leq \sum_{J \text{ not adjacent to } I}|A \cap J|.
\end{equation}
 Set $A_1 = A \cap I,$ and $$A_2 = \bigcup_{J \text{ not adjacent to } I} A \cap J.$$ The desired separation is immediate and the non-concentration conditions on $A_1$ and $A_2$ are readily checked.
\end{proof}
    Let $A \subset \R^d.$ Let $A_\delta$ denote the $\delta/2$-neighbourhood of $A.$
\begin{definition}\label{def.uniformmeasure}
Let $A \subset \R^d$ be finite and $\delta$-separated. Call $\calL^1_{|A_\delta}$ the \textit{uniform measure} on $A_\delta.$ If $X$ is a random variable which outputs values from $A_\delta,$ distributed by the uniform measure, then we say that $X$ is distributed \textit{uniformly}.
\end{definition}
An important property of these measures is that if $A$ is $(\delta,\sigma,C$)-set then the uniform measure on $A$ will be $(\sigma,2C)$-Frostman.
\begin{lemma}\label{measure}
Let $0 < \delta, \sigma < 1.$ Let $A \subset \R$ be a $(\delta,\sigma,C)$-set. Let $\mu$ be the uniform measure on $A_\delta $ Then $\mu$ is $(\sigma,2C)$-Frostman. 
\begin{proof}
Recall that by setting $r = \delta$ into the non-concentration condition imposed of $A$ we must have that $|A| \geq \delta^{-\sigma}/C.$
Let $0 < r < \delta.$ Then
\begin{displaymath}
\mu(B(x,r)) \leq \frac{2r}{\delta|A|} \leq 2C\delta^{  \sigma-1}r \leq 2Cr^{\sigma}.
\end{displaymath}
Now let $\delta \leq r \leq 1.$ Then
\begin{displaymath}
\mu(B(x,r)) = \frac{\calL_1(A_\delta \cap B(x,r))}{\delta|A|} \leq \frac{2C\delta r^\sigma|A|}{\delta|A|} = 2Cr^\sigma. 
\end{displaymath}
\end{proof}
\end{lemma}
\subsection{Radial projections}
    The \textit{radial projection centred at} $x \in \R^2$ is the map $\pi_x:\R^2\setminus\{x\} \rightarrow S^{1}$ defined by
    $$\pi_x(y) := \frac{y-x}{|y-x|}.$$
For a point $x \in X$ and a set $Y \subset \R^2$ the image $\pi_x(Y\setminus \{x\})\subset S^{1}$ gives all the unit vectors in $\R^2$ one can define from line segments between $x$ and $y$ for $y \in Y.$

\begin{definition}
    For two measures $\mu$ and $\nu$ on $\R^2$ with $\dist(\spt \mu,\spt \nu) > 0$  define the \textit{quotient measure} of $\mu$ and $\nu$ by
    \begin{displaymath}
\rho_{\mu,\nu}(A) := \int\int 1_A\bigg(\frac{y_2-x_2}{y_1-x_1}\bigg)d\mu(x_1,x_2)d\nu(y_1,y_2)
\end{displaymath} 
for all measurable $A \subset \R.$
\end{definition}

 \indent We will abbreviate $\rho_{\mu,\nu}$ to $\rho$ when it is clear what $\mu$ and $\nu$ are from context. We view $\rho$ as the measure on all the gradients one can obtain from line segments generated by pairs of points in $\spt \mu \times \spt \nu.$  We need the below result as a basis to get our strong bounds for the discretised ring theorem (Theorem \ref{sumprod}). 
\begin{prop}
    \label{frostslope}
 Let $C, \epsilon >0,$ let $0 < s \leq 1,$ and let $0 \leq t < \min\{s,1\}.$ There exists $K = K(C,\epsilon,s,t) >0$ so that the following holds.

 \indent Let $\mu_1,\mu_2,\nu_1,\nu_2$ be $(s,C)$-Frostman measures supported inside $[-10,10].$ Set $\mu := \mu_1 \times \mu_2,$ and $\nu := \nu_1 \times \nu_2.$ Suppose that $\dist(\spt\mu, \spt\nu) \geq 1/C.$ Let $\rho$ be the quotient measure of $\mu$ and $\nu.$ Then there exists a $G \subset \spt \rho$ with $$\rho(G) \geq 1-\epsilon$$ such that $\rho_{|G}$ is $(t,K)$-Frostman.
\begin{proof}
Set $X := \spt \mu$ and $Y := \spt \nu.$ Apply \cite[Corollary 2.19, Theorem 3.20]{orp} to find $K > 0$ and $E \subset X$ with $\mu(E) \geq 1 - \epsilon$ so that for all $x \in E$ there exists $F_x \subset Y,$ with $\nu(F_x) \geq 1-\epsilon,$ so that $\pi_x \nu_{|F_x}$ is $(t,K)$-Frostman. Applying $\tan: S^1 \rightarrow \R$ and noting that for all balls $B$ of radius $r,$ $\tan^{-1}(B + t)$ is contained in at most $\sim 1$ balls of radius $r,$ it follows that $\rho_{\delta_x, \nu_{|F_x}}$ is $(t,O(K))$-Frostman. Integrating on $x$ with $\mu_{|E}$ and writing $G := \{(x,y): x \in E, y \in F_x\},$ we find that $\rho_{\mu,\nu | G}$ is $(t,O(K))$-Frostman with $(\mu \times \nu) \geq (1-\epsilon)^2 \geq 1- 2\epsilon.$ Replacing $\epsilon$ with $\epsilon/2$ give the result.
\end{proof}
\end{prop}
\subsection{Discretised information inequalities}\label{section.discretised.shannon}
We first recall some basics of Shannon entropy. Let $X$ be a random variable taking values in a finite set $G.$ The \textit{Shannon entropy} is defined by
\begin{equation}
    \E(X) := \sum_{x \in G}\pr(X = x)\log \pr(X=x)^{-1},
\end{equation}
where we interpret $0\log 0 := 0.$ We have the the chain rule: for two random variables $X,Y$ we have
\begin{equation}\label{eq.chainrule}
    \E(X,Y) = \E(X) + \E(Y|X),
\end{equation}
where $Y|X$ denotes the random variable $Y$ conditioned on $X.$ We also have submodularity:
\begin{thm}\label{thm.submod}
    Let $X,Y,Z,W$ be random variables and suppose that $X$ determines $Z$ and $Y$ determines $Z,$ further suppose that $(X,Y)$ determines $W.$ Then
    \begin{equation}
        \E(Z) + \E(W) \leq \E(X) + \E(Y).
    \end{equation}
\end{thm}
\indent We now consider a discretised variant: For a sample space $\Om \subset \R^d$ the set $\calD_\delta(\Om)$ denotes the intervals of
$$\{ [\delta n_1, \delta (n_1+1))\times \cdots \times [\delta n_d, \delta (n_d+1)) : (n_1,\dots,n_d) \in \Z^d\}$$
 which intersect $\Om.$
 \begin{definition}
  Let $X$ be a random variable taking values in $\Om.$  We define the \textit{Shannon entropy} with respect to $\D_\delta(\Omega)$ by
\begin{displaymath}
    \Ed(X) = \sum_{I \in \calD_\delta(\Om)}\pr(X \in I)\log\pr(X \in I)^{-1}.
\end{displaymath}
\end{definition}
We interpret $0\log 0^{-1} = 0.$ Sometimes we may write $\E_\delta(\mu)$ when we want to emphasise the underlying measure. Shannon entropy of random variables in Euclidean space with respect to a partition have been considered in many works. See, for example, \cite[Section 5.4]{fal}. We note that for two random variables $X,Y$ we have
$$\Ed(X,Y) \leq \Ed(X) + \Ed(Y),$$
with equality if $X$ and $Y$ are independent. For two i.i.d.~random variables $X,Y$ we have
$$\Ed(X,Y) = \Ed(X) + \Ed(Y) = 2\Ed(X).$$
We denote by $\spt X$ the support of $\mu,$ $\spt \mu,$ where $\mu$ is the underlying measure. We say that $X$ is $s$\textit{-Frostman} if the underlying measure $\mu$ is $s$-Frostman.

\indent We recall four well known facts that we shall need. The first is a straightforward application of Jensen's inequality: If $f:A \subset \R \rightarrow \R$ is concave, and $p_1,\dots,p_n$ is a probability vector, then 
$$\sum_{i=1}^n p_i f(x_i) \leq f\bigg(\sum_{i=1}^n p_ix_i\bigg)$$
for all $x_1,\dots,x_n \in \R.$
\begin{lemma}[Upper and lower bounds]\label{lem.maxent}
    Let $X$ be a compactly supported random variable on $\R^d.$ Then 
    $$0 \leq \Ed(X) \leq \log N_\delta(\spt X) + O(1)$$
    for all $\delta > 0.$
    \begin{proof}
The lower bound follows since $x \log 1/x \geq 0$ for $0 <  x \leq 1.$ For the upper bound, we use that the function $f(x) = \log(x)$ is concave. Therefore, applying Jensen's inequality as above, we have
\begin{align*}
    \Ed(X) &= \sum_{I\in \D_\delta(\Omega)} \pr(X \in I) \log \pr(X\in I)^{-1}\\
    &\leq \log \sum_{I\in \D_\delta(\Omega)} 1\\
    &= \log N_\delta(\spt X) + O(1),
\end{align*}
as required. 
    \end{proof}
\end{lemma}
\indent The second is continuity. 
\begin{lemma}[Continuity]\label{lem.contofent}
    Suppose that $X$ is a random variable on a compact set $A \subset \R^d$ and let $C > 1$, $ \delta > 0.$ Then
    $$\Ed(X) \leq \E_{C\delta}(X) +  O(1),$$
    with the implicit constant depending on $C$ and $d$ only.
    \begin{proof}
        Let $X_1$ be the random variable on $\calD_\delta(A)$ which outputs the $I \in \D_\delta(A)$ for which $X \in I;$ let $X_2$ be the random variable on $\calD_{C\delta}(A)$ which outputs the $J \in \D_{C\delta(A)}$ for which $X \in J.$ We have by the chain rule \eqref{eq.chainrule}, 
        $$\E(X_1,X_2) = \E(X_2) + \E(X_1 | X_2),$$
        which leads us to 
        $$\Ed(X) \leq \E_{C\delta}(X) + \E(X_1 | X_2).$$
        Finally, for each $J \in \D_{C\delta}(A)$ the random variable $(X_1 | X_2 = J)$ has a sample space of size at most $\sim C^d,$ and so for each $J$ we have
        $$\E(X_1 | X_2 = J) \leq O(d\log C),$$
        and so taking expectation gives us
        $$\E(X_1 |X_2) \leq \log O(d \log C),$$
        and the result follows.
    \end{proof}
\end{lemma}
\indent The third is stability under large restrictions.
    \begin{lemma}[Restriction]\label{rest}
Let $\mu$ be a probability measure supported on a compact $A \subset \R.$ Let $\epsilon > 0$ and let $A' \subset A$ be such that $\mu(A') \geq 1-\epsilon.$ Then 
\begin{displaymath}
(1-\epsilon)\Ed(\mu_{|A'}) \leq \Ed(\mu).
\end{displaymath} 
\begin{proof}
   We may write $\mu$ as the convex combination:
    \begin{equation}\label{eq.entconvex}
    \mu = \mu(A')\mu_{|A'} + \mu(A\setminus A')\mu_{|A\setminus A'}.
    \end{equation}
    Since $\E_\delta$ is concave, using \eqref{eq.entconvex}, we have
    \begin{equation}
        \mu(A')\Ed(\mu_{|A'}) + \mu(A\setminus A') \Ed(\mu_{A\setminus A'}) \leq \Ed(\mu). 
    \end{equation}
    Applying the assumptions and the non-negativity of entropy we arrive at the required result.
\end{proof}
\end{lemma}
\indent The fourth gives a lower bound for the Shannon entropy of Frostman measures.
\begin{lemma}[Frostman bound]\label{frost}
Suppose that $\mu$ is $(s,C)$-Frostman on $\R^d.$ Then
\begin{displaymath}
\Ed(\mu) \geq s\log \delta^{-1} -\log C - O(1).
\end{displaymath}
\begin{proof}
We have 
\begin{align*}
\Ed(\mu) &= -\sum_{I \in \calD_\delta(\spt \mu)}\mu(I)\log\mu(I)\\
&\geq -\sum_{I \in \calD_\delta(\spt \mu)}\mu(I)\log(C\delta^s) - O(1)\\
&=  s\log \delta^{-1} - \log C - O(1).
\end{align*}
\end{proof}
\end{lemma}
\indent We require the following discretised submodular inequality.
\begin{lemma}[Discretised submodular inequality]\label{lem.dissubmod}
    Let $X,Y,Z,W$ be random variables taking values in  compact subsets of $\R^{k},\R^{l},\R^{m}, \R^n$ respectively. Fix $C > 1,$ $\delta > 0.$ Suppose each of the following:
    \begin{enumerate}
        \item If we know that the outcome of $X$ lies in $I \in \D_\delta(\R^k),$ then we are able to determine a choice of $2^m$ such $J \in \D_{C\delta}(\R^m)$ which the outcome of $Z$ will lie in;
        \item If we know that the outcome of $Y$ lies in $I \in \D_\delta(\R^l),$ then we are able to determine a choice of $2^m$ such $J \in \D_{C\delta}(\R^m)$ which the outcome of $Z$ will lie in;
        \item If we know the outcome of $X$ lies in $I \in \D_\delta(\R^k),$ and the outcome of $Y$ lies in $I' \in \D_\delta(\R^l),$ then we are able to determine a choice of $2^n$ such $J \in \D_{C\delta}(\R^n)$ which the outcome of $W$ will lie in.
        \end{enumerate}
        Then, $$\Ed(Z) + \Ed(W) \leq \Ed(X) + \Ed(Y) + O(1),$$
        where the implicit constant depends on $C,k,l,m,n$ only.
\begin{proof}
    Define the discrete random variables $X',Y'$ on the sample space $\D_\delta(\R^k), \D_\delta(\R^l)$ which output the $I \in \D_\delta(\R^k), J \in \D_\delta(\R^l)$ which the outputs of $X,Y$ lie in, respectively. Similarly, define the discrete random variables $Z',W'$ on the sample space $\D_{C\delta}(\R^m)$,$\D_{C\delta}(\R^n)$ respectively  which output the $I \in \D_{C\delta}(\R^m),J \in \D_{C\delta}(\R^n)$ which the outputs of $Z,W$ lie in, respectively. It is clear that
    $$\E(X') = \Ed(X), \qquad \E(Y') = \Ed(Y),$$
    and
    $$\E(Z') = \E_{C\delta}(Z), \qquad \E(W') = \E_{C\delta}(W).$$
    
   \indent By submodularity (Theorem \ref{thm.submod})
    $$\E(X',Y',Z') + \E(Z') \leq \E(X',Z') + \E(Y',Z').$$
     By construction, $X',$ determines $2^m$ potential choices of $Z',$ as does $Y',$ and $(X',Y')$ determines $2^n$ potential choices of $W'.$ Therefore using the above and the chain rule we obtain
\begin{align}
    \E(Z') + \E(W') &\leq \E(Z') + \E(W' |X',Y') + \E(X',Y') \\
    &\leq \E(Z') + \E(X',Y') + n \\
    &\leq \E(X',Y',Z') + \E(Z') + n \\ 
    &\leq \E(X',Z') + \E(Y',Z') + n \\
    &= \E(X') + \E(Y') + \E(Z'|X') + \E(Z'|X') +n\\
    &\leq \E(X') + \E(Y') + 2m + n.
\end{align}
    Using the above identifications gives us,
    $$\E_{C\delta}(Z) + \E_{C\delta}(W) \leq \Ed(X) + \Ed(Y) +2m + n$$
    Finally by the continuity of entropy (Lemma \ref{lem.contofent}) we have the result required.
\end{proof}
\end{lemma}
\indent An application of this is.
\begin{prop}\label{slopes}
Let $C > 1.$ Let $X,Y,Z,W,X'$ be random variables which take values in $[1,2].$ Let $Z'$ and $W'$ be random variables which take values in  $A_1,A_2 \subset [1,2]$ respectively, where $\dist(A_1,A_2) > 1/C.$  Then the following inequalities hold:
\begin{enumerate}
\item \begin{align*}
\Ed\bigg(\frac{X-Y}{Z'-W'}\bigg) &+ \Ed(X,X',Y,Z', W') \\
\leq \Ed(X-Y, Z'-W') &+ \Ed(X X',Y X',Z' X', W' X') + O(1). 
\end{align*}
If $X,Y,X'$ are i.i.d.~and $X,Y,X', Z',W'$ are independent then
\begin{align*}
 \Ed\bigg(\frac{X-Y}{Z'-W'}\bigg) &+ 3\Ed(X) + \Ed(Z') + \Ed(W')\\
 \leq 2\Ed(X Y) &+ \Ed(X  Z') + \Ed(X  W')+ \Ed(X-Y) + \Ed(Z'-W') + O(1).
\end{align*}
\item \begin{align*}
 \Ed\bigg(\frac{X+Y}{Z+W}\bigg) &+ \Ed(X,X',Y, Z, W) \\
 \leq \Ed(X+Y, Z+W) &+ \Ed(X X',Y X',Z X', W X') + O(1). 
\end{align*}
If $X,Y,Z,W,X'$ are i.i.d.~then
\begin{displaymath}
\Ed\bigg(\frac{X+Y}{Z+W}\bigg) + 5\Ed(X) \leq 4\Ed(X Y) + 2\Ed(X+Y) + O(1).
\end{displaymath}
\item \begin{align*}
 \Ed\bigg(\frac{X-Y}{Z'-W'}\bigg) &+ \Ed(X,Y, Z', W') \\
 \leq \Ed(X-Y, Z'-W') &+ \Ed(Y/X,Z'/X,W'/X) + O(1). 
\end{align*}
If $X,Y$ are i.i.d.~and $X,Y, Z',W'$ are independent then
\begin{align*}
 \Ed\bigg(\frac{X-Y}{Z'-W'}\bigg) &+ 2\Ed(X) + \Ed(Z') + \Ed(W')\\
 \leq 2\Ed(Y/X) &+ \Ed(Z'/X) + \Ed(W'/X)+ \Ed(X-Y) + \Ed(Z'-W') + O(1).
\end{align*}
\item \begin{align*}
 \Ed\bigg(\frac{X+Y}{Z+W}\bigg) &+ \Ed(X,Y, Z, W) \\
 \leq \Ed(X+Y, Z+W) &+ \Ed(Y/X,Z/X, W/X) + O(1). 
\end{align*}
If $X,Y,Z,W$ are i.i.d.~then
\begin{displaymath}
\Ed\bigg(\frac{X+Y}{Z+W}\bigg) + 4\Ed(X) \leq 3\Ed(X/Y) + 2\Ed(X+Y) + O(1).
\end{displaymath}
\end{enumerate}
\begin{proof}
    We prove the first, the rest are similar.
    
    \indent Suppose we know that $(X-Y,Z'-W') \in I \times J,$ where $I,J$ intervals of length $\delta.$ Since $\dist(Z',W') > 1/C$ we see that $\tfrac{X-Y}{Z'-W'}$ lies in an interval of length $2C\delta.$

    \indent Suppose we know that $(XX',YX',Z'X',W'X') \in I_1 \times I_2 \times I_3 \times I_4,$ each an interval of length $\delta.$ Then $$\tfrac{XX' - YX'}{Z'X' - W'X'} = \tfrac{X-Y}{Z'-W'}$$ lies in an interval of length $2C\delta.$

    \indent Now suppose that we know both of the facts stipulated at the beginning on the previous two paragraphs. Then $$X' = \tfrac{Z'X' - W'X'}{Z'-W'}$$ will lie in an interval of length $C^2\delta.$ Then each $X,Y,Z',W'$ will each lie in a (separate) interval of length $C^3\delta.$ The result then follows from Lemma \ref{lem.dissubmod}. 
\end{proof}
\end{prop}
\section{Proofs of Theorem \ref{sumprod}}
The below is a key lemma from which our results will follow easily. 
\begin{lemma}\label{lemma.keylemma}
Let $C,K,\epsilon > 0$ and let $0 < t < \min\{2\sigma,1\}.$ Let $\mu, \mu_1, \mu_2$ be $(\sigma,C)$-Frostman supported on $[1,2],$ where $\dist(\spt \mu_1, \spt \mu_2) \geq 1/K.$ Write $A = \spt \mu,$ and suppose that $\spt \mu_1, \spt \mu_2 \subset A.$ Then
\begin{align}
    (t(1-\epsilon) + 5\sigma )\log \delta^{-1} &\leq \log N_\delta((A-A)^2\times (AA)^4) + O(1);\\
        (t(1-\epsilon) + 5\sigma )\log \delta^{-1} &\leq \log N_\delta((A+A)^2\times (AA)^4)  + O(1);\\
            (t(1-\epsilon) + 4\sigma )\log \delta^{-1} &\leq \log N_\delta((A-A)^2\times (A/A)^3) + O(1);\label{eq.hard} \\
                (t(1-\epsilon) + 4\sigma )\log \delta^{-1} &\leq \log N_\delta((A+A)^2\times (A/A)^3) + O(1).
\end{align}
The implicit constants depend on $C,K,\epsilon,\sigma,t.$
\begin{proof}
     We prove \eqref{eq.hard}, the rest are similar. Let $X,Y,$ be i.i.d.~random variables distributed by $\mu,$ let $Z$ be a random variable distributed by $\mu_1,$ and let $W$ be a random variable distributed by $\mu_2.$ Since each of these random variables is $\sigma$-Frostman, it follows from Lemma \ref{frost} that their entropies at scale $\delta$ are at least $\sigma \log \delta^{-1} - O(1).$ Applying Proposition \ref{slopes}, along with this fact, we obtain
\begin{align*}
            \Ed\bigg(\frac{X-Y}{Z-W}\bigg)  + 4\sigma \log \delta^{-1} &\leq \Ed(X-Y,Z-W) + \Ed(Y/X,Z/X,W/X) + O(1).
\end{align*}
 Applying Proposition \ref{frostslope} we see that the random variable
$\tfrac{X-Y}{Z-W}$
is $t$-Frostman when conditioned to a subset of measure at least $1-\epsilon.$ Applying Lemma \ref{rest} and then Lemma \ref{frost}, along with independence, the trivial upper bounds for entropy (Lemma \ref{lem.maxent}), and the fact that $\spt \mu_1, \spt \mu_2 \subset A,$ we obtain the required inequalities.
\end{proof}
\end{lemma}
\begin{proof}[Proof of Theorem \ref{sumprod}]
We prove \eqref{eq.hard2} with $-,$ when $0 < \sigma \leq 1/2,$ the rest are similar. Let $\epsilon > 0$ and $0 < t < 2\sigma.$ Let $\mu$ be the uniform measure on $A_\delta.$
 We apply Lemma \ref{lem.twosubsets} to find a $C' \sim C$ and $A_1, A_2 \subset A,$ both $(\delta,\sigma,C')$ sets with $\dist(A_1,A_2) \gtrsim C'^{-1}.$ Let $\mu_1$ be the uniform measure on $(A_1)_\delta,$  and let $\mu_2$  be the uniform measure on $(A_2)_\delta,$  Since these measures are $(\sigma,2C')$-Frostman. Applying Lemma \ref{lemma.keylemma} and taking exponents gives us
\begin{align*}
     \delta^{-(t(1-\epsilon) + 4\sigma )} &\lesssim N_\delta(A-A)^2N_\delta(A/A)^3;
\end{align*}
Since $0 < \epsilon < 1$ and $0 < t < 2\sigma$ are arbitrary, the result follows. 
\end{proof}

\section{Examples}\label{section.examples}
We note that for a set $A \subset \R$ we cannot guarantee that 
\begin{equation}
    \max\{\dimh(A+A), \dimh(AA)\} > \dimh A,
\end{equation}
when $0 < \dimh A < 1.$
Fix $0 < s < 1.$ Let $0 < \alpha < 1$ and $\beta > 1.$ For each $N \in \N$ set 
\begin{equation}
    A_N := \{\alpha,2\alpha,\ldots,N\alpha\}
\end{equation}
and 
\begin{equation}
    G_N := \{\beta,\beta^2,\ldots , \beta^{N}\}.
\end{equation}
Consider the collections of maps $$\{cx + i\}_{i \in A_N}, \{cx + j\}_{j \in G_N},$$ where $0 < c <1.$
Associate a word $\om = (\om_1,\om_2,\ldots) \in (A_N \cup G_N)^\infty,$ to a point $x_\om,$ by the relation
\begin{equation}
x_\om = \sum_{i=1}^\infty c^i\om_i.    
\end{equation}
For each $n \in \N$ let $k_n := 10^n.$ We define an admissible word $\om \in (A_N \cup G_N)^\infty$ as follows: if $k_{2n} < k < k_{2n+1}$ then $\om_k \in A_N,$ otherwise $\om_k \in G_N.$ Call the collection of admissible words $\Om.$ Choose $c,\alpha,\beta,N$ so that $s = \log N/\log (1/c),$ and so that the cylinder sets for each $k \in \N$ are disjoint. Set
\begin{equation}
    A := \{x_\om: \om \in \Om\}.
\end{equation}
\begin{prop}
    We have $$s = \dimh A = \lbox A = \ubox A,$$ and $$\lbox (A+A) = \lbox AA = \dimh A.$$
\end{prop}
Therefore, is is not possible to expect even a gain for lower-box dimension. 

\indent Secondly, we show that \cite[Corollary 5.8]{furen} is sharp when $2/3 < \sigma \leq 1.$
\begin{prop}\label{prop.upperbound}
For all $0 < \sigma < 1$ there exists a constant $C = C(\sigma) > 0$ so that for all $\delta >0$ we may find a $(\delta,\sigma, C)$-set $B$ with $|B| \sim \delta^{-\sigma}$ and 
$$N_\delta(B+B) + N_\delta(BB) \lessapprox \delta^{-(1+\sigma)/2}.$$
\end{prop}
For $N \in \N, 0 < \alpha < 1, \beta >1,$ let $A_N$ and $G_N$ be as before.
Consider the collections of maps as before,
and let $A$ and $G$ be their respective attractors. Let $\epsilon > 0$ and $1/2 + \epsilon \leq \tau + \epsilon < 1.$  Choose $c$ and $N$ so that $\tau + \epsilon = \log N/\log (1/c),$ and so that $A$ and $G$ satisfy the open set condition. Then $\dimh A = \dimh G = \tau + \epsilon.$ By replacing $G$ with a randomly translated copy if necessary (see \cite[Theorem 8,1]{falbk}), we have that $A \cap G$ has Hausdorff dimension at most $2\tau + 2\epsilon -1.$ Set $\sigma = 2\tau -1.$ By \cite[Lemma 5.3]{bflm}, we may find $B \subset A \cap G$ which is a $(\delta,\sigma,C)$-set, for some $C \sim_{c,N} 1.$ We then have 
\begin{equation}
    N_\delta(B+B) + N_\delta(BB) \leq N_\delta(A+A) + N_\delta(GG) \leq 100\delta^{-\tau - \epsilon} = 100\delta^{-\frac{1+\sigma}{2}-\epsilon}.
\end{equation}
\bibliographystyle{alpha}
\bibliography{references}
\end{document}